\documentclass{amsproc}
\usepackage{amssymb}
\newtheorem{theorem}{Theorem}[section]

\theoremstyle{definition}
\newtheorem{definition}[theorem]{Definition}

\theoremstyle{remark}

\numberwithin{equation}{section}

\newcommand{\ot}{\otimes}
\newcommand{\recht}{\rightarrow}

\newcommand{\de}{\Delta}
\newcommand{\io}{\iota}
\newcommand{\eps}{\varepsilon}
\newcommand{\vfi}{\varphi}
\newcommand{\om}{\omega}
\newcommand{\cM}{\mathcal{M}}
\newcommand{\cL}{\mathcal{L}}
\newcommand{\cG}{\mathcal{G}}
\newcommand{\R}{\mathbb{R}}
\newcommand{\mul}{\operatorname{m}}
\newcommand{\la}{\lambda}
\newcommand{\al}{\alpha}
\newcommand{\be}{\beta}
\newcommand{\kom}{\; , \;}
\begin{document}

\title[Examples of locally compact quantum groups]{Examples of locally compact quantum groups through the bicrossed product construction}

\author{Stefaan Vaes}
\address{Dept. of Mathematics, Katholieke Universiteit Leuven, Celestijnenlaan 200~B, B--3001 Heverlee, Belgium}

\email{Stefaan.Vaes@wis.kuleuven.ac.be}

\thanks{The author is a Research Assistant of the Fund for Scientific Research Flanders --
Belgium (F.W.O.)}

\date{June 2000}

\maketitle

Building on work of Kac \& Vainerman \cite{VK}, Enock \& Schwartz \cite{ESbook}, Baaj \& Skandalis \cite{BS}, Woronowicz \cite{Wor5} and Van
Daele \cite{VD}, locally compact quantum groups were introduced by J.~Kustermans and the author in
\cite{KV2} (see also \cite{KV1} and \cite{KVPnas}).

The main aim of this paper is to introduce some examples of non-compact locally compact quantum
groups to a non-specialized audience. The major importance of these examples is their
simplicity. Other examples as the quantum $E(2)$ group of Woronowicz (see \cite{Wor2}) are much more
difficult to construct.
We will make use of the bicrossed product construction to obtain our examples. This
construction dates back to Kac \cite{Kac}, Takeuchi \cite{Tak} and Majid \cite{Maj} and \cite{Maj2}. It was also considered by Baaj \& Skandalis
in \cite{BS2}.

At the same time we will try to explain what is a locally compact quantum group without going
into details and we will sketch why to study locally compact quantum groups.
It is our hope that this paper can contribute to make locally compact quantum groups or the
operator algebra approach to quantum groups in general, more attractive to a broader audience.

\subsection*{Acknowledgment}
I would like to express my gratitude towards Leonid Vainerman for
numerous discussions on the bicrossed product construction and the examples coming out of it.

\section{Why locally compact quantum groups?}
In this section we want to give three motivations for the study of locally
compact quantum groups.
\subsection{The duality of groups}
Let $G$ be an Abelian locally compact group. Then the set of continuous characters on
$G$ can be given a locally compact group structure $\hat{G}$. The famous Pontryagin-Van Kampen
theorem says that the bidual group is canonically isomorphic with $G$. Trying to generalize
this duality to non-Abelian groups one runs into serious difficulties.

When $G$ is a compact group, one can build a difficult algebraic structure starting with the
irreducible representations of $G$. One can think e.g. of the Tannaka-Krein duality. Of course
this dual structure will not be a group anymore. When $G$ is no longer compact, things become
even more difficult, because now some irreducible representations can be infinite dimensional.
In this situation algebraic methods will not suffice any more.

So, it is an interesting problem to look for a self-dual category containing all locally
compact groups. This problem was solved independently by Kac \& Vainerman (see \cite{VK}) and
Enock \& Schwartz (see \cite{ESbook}) in the seventies. The object they defined is nowadays
called a Kac algebra and the theory of Kac algebras can be considered as the first topological
theory of quantum groups.

With the development of compact quantum groups by Woronowicz, it became clear that the category
of Kac algebras was too small to include all quantum groups. For instance the quantum $SU(2)$
group cannot be obtained as a Kac algebra. So it became a major challenge to unify compact
quantum groups and Kac algebras in one theory of locally compact quantum groups.
\subsection{Quantum symmetries}
In the important papers \cite{EN} and \cite{E} Enock and Nest have shown that locally compact
quantum groups appear naturally as quantum symmetries of quantum spaces, in the study of
subfactors.

When a finite group $G$ acts outerly on a factor $M$ one can consider the fixed point algebra
$M^G$ as a subfactor of $M$. The subfactor $M^G \subset M$ is irreducible and of depth~2. Then
it is a natural question whether all irreducible, depth~2 subfactors are of this form. The
answer is no: one has to allow locally compact quantum group actions. The theorem of Enock and
Nest says that every irreducible, depth~2 inclusion of factors (possibly with infinite index),
satisfying a regularity condition, is of the form $M^\cG \subset M$ where $\cG$ is a locally
compact quantum group acting outerly on a factor $M$.
\subsection{Topological version of Hopf algebras}
In some parts of mathematical physics Hopf algebras play the role of quantum symmetries.
Sometimes Hopf algebras are simply called quantum groups. In our opinion this is not completely
justifiable.
\begin{itemize}
\item The theory of Hopf algebras does not allow an analogue of the Haar measure. Taking into
account the major role played by the Haar measure in harmonic analysis, we consider this as a
major drawback.
\item One cannot study infinite dimensional corepresentations of Hopf algebras.
The matrix coefficients of these infinite dimensional corepresentations will not be in the Hopf
algebra, but will rather be continuous functions of elements in the algebra.
\item There is no notion of positivity in a Hopf algebra.
\item We have thrown away all topology.
\end{itemize}
In the theory of locally compact quantum groups we will present in a sense a topological
definition of a Hopf algebra. The underlying algebra will be replaced by a C$^*$-algebra or von
Neumann algebra. Instead of looking at polynomials we look at continuous functions that vanish
at infinity, or at bounded measurable functions. Then the objections made above can be
handled.

\section{The definition of a locally compact quantum group}
Before I state the definition of a locally compact quantum group as it was given by
J.~Kustermans and the author in \cite{KV2} and \cite{KV3}, I want to motivate where it comes from. So, let $G$
be a locally compact group. Then consider the von Neumann algebra $M=L^\infty(G)$. The group
law of $G$ can be translated to a homomorphism
$$\de:M \recht M \ot M : (\de f)(g,h) = f(gh)$$
where $\ot$ denotes the von Neumann algebraic tensor product, and where we identified $M \ot M$
with $L^\infty(G \times G)$. The associativity of the group law can be translated to the
formula
$$(\de \ot \io)\de = (\io \ot \de) \de$$
called coassociativity, where $\io$ denotes the identity map. Now the philosophy will be to replace the commutative von Neumann
algebra $L^\infty(G)$ by an arbitrary von Neumann algebra $M$. Hence the starting point is a
pair $(M,\de)$ of a von Neumann algebra $M$ and a homomorphism
$$\de: M \recht M \ot M \quad\text{satisfying}\quad (\de \ot \io)\de = (\io \ot \de)\de.$$
Of course, this is not enough to call the pair $(M,\de)$ a quantum group. We have only
translated an associative group law, and so we only have a quantum semigroup. We have to add
extra axioms to obtain quantum groups.

The obvious guess for other axioms is of course to add the existence of antipode and counit.
This has been tried by many mathematicians for many years, and it failed. There are several
obstructions in defining the right notion of antipode and counit in this operator algebra
setting.
\begin{itemize}
\item The antipode and counit can be unbounded, not everywhere defined.
\item The multiplication map $\mul : M \ot M \recht M$ can be unbounded for non-commutative
$M$.
\end{itemize}
So, the usual Hopf algebra axiom characterizing the antipode $S$:
$$\mul (S \ot \io) \de(a) = \eps(a) 1 = \mul (\io \ot S) \de(a)$$
cannot be given a meaning. Therefore one has to look for alternative axioms. In \cite{KV2} and
\cite{KV3} J.~Kustermans and the author discovered that it is enough to assume the existence of
invariant Haar measures. So we have given the following definition.
\begin{definition} \label{defi}
A pair $(M,\de)$ is called a locally compact quantum group, when $M$ is a von Neumann algebra
and $\de$ is a homomorphism
$$\de: M \recht M \ot M \quad\text{satisfying}\quad (\de \ot \io)\de = (\io \ot \de)\de.$$
Further we assume the existence of normal semifinite faithful weights $\vfi$ and $\psi$ on $M$
satisfying
\begin{align}
&\vfi \bigl( (\om \ot \io) \de(a) \bigr) = \om(1) \, \vfi(a) \quad\quad a \in \cM_\vfi^+,\om \in M_*^+ \label{links} \\
&\psi \bigl( (\io \ot \om) \de(a) \bigr) = \om(1) \, \psi(a) \quad\quad a \in \cM_\psi^+,\om \in
M_*^+. \label{rechts}
\end{align}
\end{definition}
Let us look somewhat closer at
this definition. We denoted with $\cM_\vfi^+$ the set of elements $a \in M^+$ satisfying
$\vfi(a) < +\infty$. Then formulas~\ref{links} and \ref{rechts} express the
left and right invariance of the weights $\vfi$ and $\psi$. Hence $\vfi$ and $\psi$ are the
analogues of the Haar measure on a locally compact group. We will sometimes call them Haar
weights. This will become clear in the examples in the next section.

Also it should be stressed that the definition given above is stated in the von Neumann
algebra language (as it was done in \cite{KV3}). In \cite{KV2} we have given a definition in
the C$^*$-algebra language, and from a philosophical point of view C$^*$-algebras should be
preferred. Indeed, the C$^*$-algebra reflects the locally compact topology, while the von
Neumann algebra reflects rather the measurable structure. Nevertheless it should be emphasized
that both approaches are equivalent: there is a canonical bijective correspondence between
C$^*$-algebraic and von Neumann algebraic quantum groups. This can be considered as a quantum
version of the classical Weil theorem, saying that every group with left invariant measure
class can be given a locally compact topology making it a locally compact group.

Before turning to some examples of locally compact quantum groups, let us explain what kind of
theory can be developed from this definition. In the definition we assumed the existence of
Haar weights, but we can prove their uniqueness (up to a positive scalar). We can also prove the
existence of a densely defined, possibly unbounded antipode $S$, satisfying
$$S \bigl( (\io \ot \vfi) (\de(a)(1 \ot b)) \bigr) = (\io \ot \vfi)((1 \ot a) \de(b)).$$
This antipode $S$ will have a polar decomposition $S=R \tau_{-i/2}$, where $R$ is an
anti-automorphism of $M$ and $(\tau_t)_{t \in \R}$ is a one-parameter group of automorphisms of
$M$. Then $R$ is called the unitary antipode and $\tau$ the scaling group of $(M,\de)$. We can
also construct a dual locally compact quantum group $(\hat{M},\hat{\de})$ and prove a
Pontryagin-Van Kampen duality theorem. There exists a positive number $\nu > 0$, called the
scaling constant, such that $\vfi \tau_t = \nu^{-t} \vfi$ for all $t \in \R$. Recently Woronowicz and
Van Daele have discovered the first example of a locally compact quantum group with scaling constant
different from $1$: the quantum $az+b$ group.

Of course it would be desirable to give a definition of locally compact quantum groups without
assuming the existence of invariant weights, and with their existence as a theorem. This aim
seems to be far out of reach, and not all researchers believe it can be reached. Further, in
the construction of examples the Haar measure does not cause great problems, and usually it is
just there.
\section{Examples of locally compact quantum groups}
Before giving some non-trivial examples, we will look at the examples coming from locally
compact groups. These examples mainly help to interpret locally compact quantum group results.
\subsection{Locally compact quantum groups coming from groups}
Let $G$ be a locally compact group. Define $M=L^\infty(G)$. As before we put $\de: M \recht M
\ot M$, defined by $(\de f)(g,h) = f(gh)$. One can define a normal semifinite faithful weight
$\vfi$ on $M$ by the formula
$$\vfi(f) = \int f(x) \; dx$$
where we integrate with respect to the left Haar measure on $G$. It is easy to check that
$\vfi$ is a left invariant weight in the sense of formula~\ref{links} of definition~\ref{defi}.
Analogously one can obtain a right invariant weight in the sense of formula~\ref{rechts}. Then
the pair $(M,\de)$ will be a locally compact quantum group, and it will be the most general
commutative locally compact quantum group.

Next we can put $M = \cL(G)$, the group von Neumann algebra. When $(\la_g)_{g \in G}$ is the
left regular representation of $G$ on $L^2(G)$, the von Neumann algebra $M$ has a dense
$*$-subalgebra consisting of the elements
$$\int f(g) \, \la_g \; dg$$
where $f$ is a continuous, compactly supported function. Then one can construct a normal
semifinite faithful weight $\vfi$ on $M$ such that
$$\vfi \Bigl( \int f(g) \, \la_g \; dg \Bigr) = f(e).$$
On a slightly formal level it is easy to check the invariance of $\vfi$. Let $\om \in M_*$
and let $f$ be a continuous compactly supported function. Then we have
\begin{align*}
\vfi\Bigl( (\om \ot \io) \de \Bigl( \int f(g) \lambda_g \; dg \Bigr) \Bigr)
&= \vfi \Bigl( \int f(g) \om(\lambda_g) \lambda_g \; dg \Bigr) \\
&= f(e) \; \om(\lambda_e)
= \vfi \Bigl( \int f(g) \lambda_g \; dg \Bigr) \; \om(1).
\end{align*}
So, again $(M,\de)$ is a locally compact quantum group, and it is the most general
cocommutative locally compact quantum group, i.e. satisfying $\sigma \de = \de$, where
$\sigma$ flips the two legs of $M \ot M$.

Finally we mention that $\cL(G)$ will be the dual of $L^\infty(G)$ and vice versa.
\subsection{Non group-like examples}
Several examples of non-compact locally compact quantum groups have been constructed by
Woronowicz. The most important ones are the quantum $E(2)$ group, the quantum $SL(2)$ group and
the quantum $ax+b$ and $az+b$ groups. The Haar weight on the quantum $E(2)$ group has been
constructed by Baaj and the Haar weight on the quantum $ax+b$ and $az+b$ groups is due to Van
Daele. The main problem with these examples is the difficulty of their construction.
Introducing a single example requires a lot of difficult, but beautiful mathematics.

So, the major aim of the next two paragraphs will be to introduce another, and easily constructible
example. Nevertheless this example is non-trivial, which roughly means that it is not a Kac
algebra (in particular it is not group-like) and that it is non-compact and non-discrete.
\subsection{Locally compact quantum groups through bicrossed products}
The starting point for this construction is a matched pairs of locally compact groups. Matched
pairs of finite groups were defined by G.~Kac in \cite{Kac}, and used to obtain the historical
Kac-Paljutkin dimension~8 example of a Kac algebra. These matched pairs were rediscovered by
Takeuchi in \cite{Tak} and Majid gave a definition for locally compact groups in \cite{Maj2}.
This construction has also been used by Baaj \& Skandalis in \cite{BS2}.

Let $G$ and $H$ be locally compact groups and suppose that the maps
\begin{align*}
& G \times H \recht H : (g,s) \mapsto \al_g(s) \\
& G \times H \recht G : (g,s) \mapsto \be_s(g)
\end{align*}
are defined nearly everywhere and measurable. Suppose further that
\begin{alignat*}{3}
\be_{st}(g) &= \be_s(\be_t(g)) \kom &\qquad \al_g(st) &= \al_{\be_t(g)}(s) \al_g(t) \kom &\qquad \text{for nearly all }\; & (s,t,g) \kom \\
\al_{gh}(s) &= \al_g(\al_h(s)) \kom &\qquad \be_s(gh) &= \be_{\al_h(s)}(g) \be_s(h) \kom &\qquad \text{for nearly all }\; & (g,h,s).
\end{alignat*}
Then we call $(G,H)$ a matched pairs of groups. In \cite{Maj2} Majid supposes that $\al_g(s)$
and $\be_s(g)$ are everywhere defined and continuous. But to obtain the interesting examples we
will discuss further, this is too restrictive.
Suppose now that such a matching is given. Then we can define
$$\al: L^\infty(H) \recht L^\infty(G \times H) : (\al f)(g,s) = f(\al_g(s)).$$
Denoting with $\de_G$ the usual comultiplication on $L^\infty(G)$ we have
$$(\io \ot \al)\al = (\de_G \ot \io) \al.$$
Hence $\al$ will be an action of $G$ on the von Neumann algebra $L^\infty(H)$. So we can define
the crossed product von Neumann algebra $M$ as the von Neumann algebra generated by
$\al(L^\infty(H))$ and $\cL(G) \ot 1$,
where $\cL(G)$ is the group von Neumann algebra as considered above. We observe that $M$ acts
naturally on the Hilbert space $L^2(G \times H)$.
Finally we can introduce a unitary $W$ on $L^2(G \times H \times G \times H)$ by
$$(W \xi)(g,s,h,t) = \xi( \be_{\al_g(s)^{-1} t}(h) g, s, h, \al_g(s)^{-1} t ).$$
Then one can prove that the formula
$$\de(z) = W^* (1_{L^2(G \times H)} \ot z) W \quad\text{for }\; z \in M$$
defines a comultiplication on $M$ turning $(M,\de)$ into a locally compact quantum group. The
left invariant weight on $(M,\de)$ can be obtained as the dual weight on the crossed product
$M$ of the left invariant weight on $L^\infty(H)$.
\subsection{Concrete example}
Let $G$ be the group $(\R,+)$ and define
$$H = \{ (a,b) \mid a > 0, b \in \R \} \quad\text{with}\quad (a,b) \, (x,y) = (ax, ay +
\frac{b}{x}).$$
Also define
$$
\alpha_x(a,b) = \begin{cases} (a+bx,b) \;\text{ if }\; a+bx > 0 \\ (-a-bx,-b) \;\text{ if }\;
a+bx < 0 \end{cases} \quad\text{and}\quad
\beta_{(a,b)}(x) = \frac{x}{a(a+bx)}.
$$
Then it is an easy exercise to check that $(G,H)$ is a matched pair of groups and the
construction above yields a locally compact quantum group $(M,\de)$. The underlying von Neumann
algebra is a crossed product and it is generated by bounded functions of the following
unbounded operators defined on $L^2(G \times H)$.
\begin{align*}
(A \xi)(x,a,b) &= \begin{cases} (a+bx) \xi(x,a,b) \;\text{ if }\; a+bx > 0 \\ -(a+bx) \xi(x,a,b) \;\text{ if }\;
a+bx < 0 \end{cases} \\
(B \xi)(x,a,b) &= \begin{cases} b \xi(x,a,b) \;\text{ if }\; a+bx > 0 \\ -b \xi(x,a,b) \;\text{ if }\;
a+bx < 0 \end{cases} \\
(C \xi) (x,a,b) &= (\partial_1 \xi)(x,a,b),
\end{align*}
where $\partial_1$ denotes the partial derivative to the first variable.
Bounded functions of the operators $A$ and $B$ generate $\al(L^\infty(H))$ and bounded functions
of the operator $C$ generate $\cL(G) \ot 1$.
Making formal computations we see that the operators $A,B$ and $C$ generate a $*$-algebra, with relations
$$A=A^* \;,\; B=B^* \;,\; C=-C^* \; , \quad [A,C]=B \;,\quad B \;\text{ central}.$$
So we get the universal enveloping algebra of the Heisenberg Lie algebra. Continuing this
formal computation we get the following comultiplication
$$
\de(A) = A \ot A \; , \quad \de(B) = A \ot B + B \ot A^{-1} \; , \quad \de(C) = C \ot A^{-2} + 1 \ot C.
$$
Hence we get a Hopf $*$-algebra, with antipode and counit given by
$$\eps(A)=1 \kom \eps(B)=\eps(C)=0 \kom S(A)=A^{-1} \kom S(B)= -B \kom S(C) = -C A^2.$$
By symmetry it is clear that also $(H,G)$ with $\be$ and $\al$ is a matched pair of groups. The
resulting locally compact quantum group will be the dual locally compact quantum group of the
previous example. Also in this case one can compute formally an underlying Hopf $*$-algebra. It
is generated by elements $A,B$ and $C$ satisfying
$$A^*=-A \kom B^*=-B \kom C^*=C \kom [A,B]=2B \kom [A,C]=2C \kom [B,C]=C^2.$$
The comultiplication is given by
$$\de(A) = A \ot 1 + 1 \ot A \kom \de(B)=B \ot 1 + 1 \ot B + A \ot C \kom \de(C)=C \ot 1 + 1
\ot C.$$
Finally antipode and counit are given by
$$\eps(A)=\eps(B)=\eps(C)=0 \kom S(A)=-A \kom S(B)=-B + AC \kom S(C)=-C.$$
In the previous example the underlying algebra was the universal enveloping algebra of the
Heisenberg Lie algebra and now the underlying coalgebra structure precisely reflects the
structure of the Heisenberg group with multiplication law
$$(a,b,c) \; \cdot \; (a',b',c') = (a+a', b+b'+ac' , c+c').$$
This is typical for the duality of quantum groups.

\bibliographystyle{amsplain}

\end{document}